\documentclass[12pt]{article}
\usepackage{amsfonts}
\usepackage{mathrsfs}
\usepackage{amssymb,amsfonts,amsmath, psfrag,eepic,colordvi,epsfig}
\numberwithin{equation}{section}
\newtheorem{theorem}{Theorem}[section]

\newtheorem{corollary}[theorem]{Corollary}

\newtheorem{conjecture}[theorem]{Conjecture}

\newtheorem{lemma}[theorem]{Lemma}

\def\pf{\noindent {\it Proof.} }
\def\qed{\hfill \rule{4pt}{7pt}}

\newcommand{\Var}{{\rm Var}}

\begin{document}
\parskip 10pt

\begin{center}

{\Large\bf  The Limiting Distribution of

the Coefficients of  the $q$-Catalan Numbers}

\vskip 6mm

William Y.C. Chen$^1$, Carol J. Wang$^2$,  and Larry X.W. Wang$^3$ \\
Center for Combinatorics, LPMC-TJKLC \\
Nankai University, Tianjin 300071, P. R. China

\vskip 3mm

$^1$chen@nankai.edu.cn, $^2$wangjian@cfc.nankai.edu.cn,
$^3$wxw@cfc.nankai.edu.cn

\end{center}
\begin{abstract}
We show that the limiting distributions of the coefficients of the
$q$-Catalan numbers and the generalized $q$-Catalan numbers are
normal. Despite the fact that these coefficients are not unimodal
for small $n$, we conjecture that for sufficiently large $n$, the
coefficients are unimodal and even log-concave except for a few
terms of the head and tail.
\end{abstract}

\noindent {\bf Keywords :} Bernoulli number, $q$-Catalan number,
unimodality, log-concavity, moment generating function.

\section{Introduction}\label{sec-int}

The main objective of this paper is to show that the limiting
distribution of the coefficients of the $q$-Catalan numbers is
normal. The Catalan numbers
\[C_n=\frac{1}{n+1}{2n\choose n}\] have many combinatorial
interpretations,  see Stanley \cite{stanbook}. The usual $q$-analog
of the Catalan numbers is given by
\begin{equation}\label{C2}
C_n(q):=\frac{1}{[n+1]}{2n\brack n},
\end{equation}
where $[n]=1+q+q^2+\cdots +q^{n-1}$, and
\[ {n\brack
k}=\frac{[n]!}{[k]! [n-k]!}.\] There are also other types of
$q$-analogs of the Catalan numbers, see, for example,
 Andrews
\cite{Andrews}, Gessel and Stanton \cite{Ges},   Krattenthaler
\cite{Kra}.

We also consider the limiting distribution of the coefficients of
the quotient of two products of $q$-numbers, which includes the
result for the $q$-Catalan numbers as a special case. We conclude
this paper with two conjectures on the unimodality and log-concavity
for almost all the coefficients of the $q$-Catalan numbers and the
generalized $q$-Catalan numbers provided that $n$ is sufficiently
large.

\allowdisplaybreaks

\section{The Limiting Distribution}\label{sec-q}

In this section, we  use the moment generating function technique to
obtain the limiting distribution of the coefficients of the
$q$-Catalan numbers. We introduce the random variable $\xi_n$
corresponding to the probability generating function
\[  \phi_n(q)=C_n(q)/C_n.\]
As far as the computations are concerned, we will not need the
following combinatorial interpretation of $C_n(q)$. However,  for
completeness and for the sake of presentation, we would mention that
$\xi_n$ reflects the distribution of the major indices of Catalan
words of length $2n$, see, for example, \cite{FH}. We  write
\[ C_n(q) =\sum m_n(k) q^k.\]

The following lemma is concerned with the expectation and variance
of $\xi_n$.

\begin{lemma}\label{ev} We have
\begin{equation} \label{s2}
E(\xi_n)=\frac{n(n-1)}{2}\quad \mbox{and} \quad
\Var(\xi_n)=\frac{n(n-1)(n+1)}{6}.
\end{equation}
\end{lemma}

\pf By the definition of $C_n(q)$, it is easy to check the following
symmetry property of $m_n(k)$:
$$m_n(k)=m_{n}(n(n-1)-k).$$ Hence
$$E(\xi_n)=\frac{n(n-1)}{2}.$$
Let
\[
F=F(q)=\prod_{i=1}^{n-1} (1+q+\cdots +q^{n+i})\quad \mbox{and} \quad
G=G(q)=\prod \limits_{i=1}^{n-1}(1+q+\cdots +q^i). \] It is easily
verified that $C_n(q)=F/G$. Since
\begin{eqnarray*}
\left. C_n(q)'' \right|_{q=1} & = & \left.
\left(\frac{F''}{G}-\frac{FG''}{G^2}-\frac{2G'F'}{G^2}
+\frac{2G'^2F}{G^3}\right)\right|_{q=1}\\[8pt]
&  = &  \frac{1}{12} n(n-1)(3n^2-n-4)\cdot C_n,
\end{eqnarray*}
we obtain
\[
 \Var(\xi_n) =  \frac{\left.
C_n(q)''\right|_{q=1}}{C_n}+E(\xi_n)-E(\xi_n)^2 =
\frac{1}{6}n(n-1)(n+1).
\]
This completes the proof. \qed

\begin{lemma}\label{asy1} When $n \rightarrow \infty$,
we have
\[
\sum_{k=2}^{\infty}
B_{2k}\frac{t^{2k}}{2k(2k)!\sigma^{2k}}\sum_{i=2}^{n}\left((n+i)^{2k}-i^{2k}\right)\rightarrow
0
\]
uniformly for $t$ from any bounded set, where $B_j$'s are the
Bernoulli numbers and $\sigma^2$ is the variance of $\xi_n$ as given
in (\ref{s2}).
\end{lemma}

\pf The second summation can be expanded as follows:
\begin{eqnarray*}
\sum_{i=2}^{n}\left((n+i)^{2k}-i^{2k}\right)  =  \sum_{i=2}^{n}
\sum_{j=1}^{2k}{2k\choose j} n^j i^{2k-j}
 =  \sum \limits_{j=1}^{2k} {2k\choose j}\left(\sum
\limits_{i=2}^{n} n^j i^{2k-j}\right).
\end{eqnarray*}
For $k>1$, the second factor in the preceding summation is bounded
by the following integral:
\[
\sum \limits_{i=2}^{n} n^j i^{2k-j}<n^j \int_{1}^{n+1}t^{2k-j}
dt=n^j\cdot \frac{(n+1)^{2k-j+1}-1}{2k-j+1}.
\]
Consequently,
\[
\sum_{i=2}^{n}\left((n+i)^{2k}-i^{2k}\right)<2^{2k}(n+1)^{2k+1}<8^{2k}n^{2k+1}.
\]
Since $\sigma^2=\frac{n^3-n}{6}>\frac{n^3}{8}$ when $n$ is
sufficiently large, we have
\[
\sigma^{-2k}\sum_{i=2}^{n}\left((n+i)^{2k}-i^{2k}\right)<64^{2k}n^{1-k}\leq
n^{-1/3}64^{2k}n^{-k/3},
\]
for large $n$ and $k>1$. Thus
\begin{eqnarray*}
\lefteqn{ \left|\sum \limits_{2\nmid k,k\geq3}B_{2k}
\frac{t^{2k}}{2k(2k)!\sigma^{2k}}
\sum_{i=2}^{n}\left((n+i)^{2k}-i^{2k}\right)\right| \qquad }\\
&<&n^{-1/3}\sum \limits_{2\nmid
k,k\geq3}|B_{2k}|\frac{t^{2k}}{2k(2k)!}64^{2k}n^{-k/3}\\
& = &n^{-1/3}\sum \limits_{2\nmid
k,k\geq3}|B_{2k}|\frac{(64tn^{-\frac{1}{6}})^{2k}}{2k(2k)!}.
\end{eqnarray*}
In view of the following asymptotic expansion of the Bernoulli
numbers
$$|B_{2n}|\sim \frac{2(2n)!}{(2\pi)^{2n}},$$
the convergent radius $R$ of the series $\sum \limits_{2\nmid
k,k\geq3}|B_{2k}|\frac{t^{2k}}{2k(2k)!}$ equals $2\pi$. Since $t$ is
from a bounded set, when $n$ is large enough, the series
\[
\sum \limits_{2\nmid
k,k\geq3}|B_{2k}|\frac{(64tn^{-\frac{1}{6}})^{2k}}{2k(2k)!}
\]
converges. Moreover, it is evident that $64tn^{-\frac{1}{6}}<1$, we
can bound the above summation by the constant
\[
M_1=\sum \limits_{2\nmid k,k\geq3}|B_{2k}|\frac{1}{2k(2k)!}. \]
Similarly, it can be deduced that
\[
\sum \limits_{2| k,k\geq
2}B_{2k}\frac{t^{2k}}{2k(2k)!\sigma^{2k}}\sum_{i=2}^{n}\left((n+i)^{2k}-i^{2k}\right)<\frac{M_2}{n^{\frac{1}{3}}},
\]
where $M_2=\sum \limits_{2\mid k,k\geq 2}B_{2k}\frac{1}{2k(2k)!}$ is
a constant. Hence
\[
\sum_{k=2}^{\infty} B_{2k}\frac{t^{2k}}{2k(2k)!\sigma^{2k}}
\sum_{i=2}^{n}\left((n+i)^{2k}-i^{2k}\right)
<\frac{M_1+M_2}{n^{1/3}},
\]
which tends to zero as $n\rightarrow \infty$. This completes the
proof. \qed

\begin{theorem}\label{main}
When $n\rightarrow \infty$, the random variable
\[
\eta_n=\frac{\xi_n-E(\xi_n)}{{\Var({\xi_n})}^{\frac{1}{2}}}
\]
has the standard normal distribution.
\end{theorem}

\pf Let $M_n(q)$ denote the moment generating function of $\xi_n$.
 Then we have $M_n(q) = \phi_n(e^q)$, see
 Sachkov \cite{Sac}. Hence
\begin{eqnarray*}
M_n(q) &= & \frac{n+1}{{2n\choose n}}
\frac{1-e^{q}}{1-e^{(n+1)q}}\cdot
\prod_{i=1}^{n}\frac{1-e^{(n+i)q}}{1-e^{iq}}\\[5pt]
& = & \prod_{i=2}^{n} \frac{i}{n+i} \cdot \prod_{i=2}^{n} \frac{1-e^{(n+i)q}}{1-e^{iq}}\\[5pt]
& = & \prod_{i=2}^{n}\frac{(1-e^{(n+i)q})/(n+i)}{(1-e^{iq})/i}\\[5pt]
& = &
\exp\left\{\frac{1}{2}\sum_{i=2}^{n}\left((n+i)q-iq\right)\right\}\prod_{i=2}^{n}
\frac{(e^{(n+i)q/2}-e^{-(n+i)q/2})/\frac{n+i}{2}}{(e^{iq/2}-e^{-iq/2})/\frac{i}{2}}\\[5pt]
& = & \exp\left\{\frac{n(n-1)q}{2}\right\} \prod \limits_{i=2}^{n}
\frac{\sinh\left((n+i)q/2\right)
/\frac{n+i}{2}}{\sinh\left(iq/2\right) /\frac{i}{2}}.
\end{eqnarray*}
Recalling the following relation  on the Bernoulli numbers
\cite{Mar}
\begin{equation}\label{bernoulli}
\ln \left(\frac{\sinh(x/2)}{x/2}\right)=\sum \limits_{k=1}^{\infty}
B_{2k} \frac{x^{2k}}{2k(2k)!},
\end{equation}
we find that
\begin{eqnarray*}
\ln M_n(q) & = & \frac{n(n-1)}{2}q+\sum \limits_{i=2}^{n} \left(\ln
\left(\frac{\sinh((n+i)q/2)}{(n+i)/2}\right)-\ln
\left(\frac{\sinh(iq/2)}{i/2}\right)\right)\\[5pt]
& = & \frac{n(n-1)}{2}q+\sum_{k=1}^{\infty}
B_{2k}\frac{q^{2k}}{2k(2k)!}\sum_{i=2}^{n}\left((n+i)^{2k}-i^{2k}\right).\\
\end{eqnarray*}
Setting $q=t/\sigma$, where $\sigma$ is the standard deviation of
$\xi_n$ as given in Theorem \ref{ev}, we are led to the expansion
\[
\ln M_n(t/\sigma)=\frac{n(n-1)t}{2\sigma}+\sum
\limits_{k=1}^{\infty} B_{2k}\frac{t^{2k}}{2k(2k)!\sigma^{2k}}\sum
\limits_{i=2}^{n}\left((n+i)^{2k}-i^{2k}\right).\] Applying Lemma
\ref{asy1},   we have, when $n\rightarrow \infty$,
\[
\sum_{k=2}^{\infty}
B_{2k}\frac{t^{2k}}{2k(2k)!\sigma^{2k}}\sum_{i=2}^{n}((n+i)^{2k}-i^{2k})\rightarrow
0
\]
uniformly for $t$ from any bounded set. Finally,
\begin{eqnarray*}
\lefteqn{\lim_{n\rightarrow \infty}M_n(t/\sigma)
\exp\left\{-\frac{n(n-1)t}{2\sigma}\right\}\quad }\\[5pt]
& = & \lim_{n\rightarrow \infty}\exp\left\{\sum_{k=1}^{\infty}
B_{2k}\frac{t^{2k}}{2k(2k)!\sigma^{2k}}\sum_{i=2}^{n}\left((n+i)^{2k}-i^{2k}\right)\right\}\\[5pt]
& = & \lim_{n\rightarrow
\infty}\exp\left\{B_2\frac{t^2}{2(2)!\sigma^2}\sum_{i=2}^{n}\left((n+i)^2-i^2\right)\right\}\\[5pt]
& = & e^{t^2/2},
\end{eqnarray*}
which coincides with the moment generating function of the standard
normal distribution. Employing  Curtiss' theorem \cite{Sac}, we
reach the conclusion that  $\eta_n$ has the standard normal
distribution when $n$ approaches infinity.\qed

\section{A General Setting}\label{sec-Gen}

 In this section, we will  determine the
 limiting distribution of the coefficients of a quotient of
 products of $q$-numbers and will give two special cases.

\begin{theorem}\label{gen}
Let $a_1, a_2, a_3, \ldots $ and $b_1, b_2, b_3, \ldots$ be two
sequences of positive numbers, and let
\[ \phi_n(x)=\sum_{k}p_n(k)x^k = \frac{(1-q^{a_1})(1-q^{a_2})\cdots
(1-q^{a_n})}{(1-q^{b_1})(1-q^{b_2})\cdots (1-q^{b_n})}.
\]
Suppose that $\xi_n$ is the random variable corresponding to the
generating function $\phi_n(x)$, that is,
\[
P(\xi_n=k)=\frac{p_n(k)}{\sum \limits_{k}p_n(k)}.
\]
Then  $\xi_n$ is normally distributed as $n\rightarrow \infty$,
 if
and only if for $k>1$
\[
\sum \limits_{k=1}^{\infty} B_{2k}\frac{t^{2k}}{2k(2k)!}\left(\sum
\limits_{i=1}^{n}(a_i^{2k}-b_i^{2k})\right)\frac{1}{\left(\sum
\limits_{i=1}^{n}(a_i^2-b_i^2)\right)^k}\rightarrow 0 \quad
\mbox{as} \quad n\rightarrow \infty.
\]
\end{theorem}

\noindent {\it Proof.} The expectation of $\xi_n$ is easy to
compute, as given below:
\[
E(\xi_n)=\phi_n(x)'_{q=1}=\frac{1}{2}\sum
\limits_{i=1}^{n}\left(a_i-b_i\right).
\]
Proceeding analogously as in the proof of Theorem \ref{ev}, we find
\begin{equation} \label{var}
\sigma^2=\Var(\xi_n)=\frac{1}{12}\sum
\limits_{i=1}^{n}\left(a_i^2-b_i^2\right).
\end{equation}
Hence,
\[
B_2\frac{t^2}{2(2)!\sigma^2} \left(\sum
\limits_{i=1}^{n}(a_i^2-b_i^2)\right)=\frac{1}{6}\cdot
\frac{t^2}{4\cdot \frac{1}{12}\left(\sum
\limits_{i=1}^{n}(a_i^2-b_i^2)\right)}\cdot \left(\sum
\limits_{i=1}^{n}(a_i^2-b_i^2)\right)=\frac{t^2}{2}.
\]
By the same procedure as in the proof of Theorem \ref{main}, we
obtain
\begin{eqnarray*}
\lefteqn{\lim \limits_{n\rightarrow \infty}M_n(t/\sigma)\exp
\left\{\frac{1}{2}\sum
\limits_{i=1}^{n}\left(a_i^{2k}-b_i^{2k}\right)\right\}}
\\ & = &
e^{t^2/2}\lim \limits_{n\rightarrow \infty} \exp \left\{\sum
\limits_{k=2}^{\infty}
B_{2k}\frac{t^{2k}}{2k(2k)!{\sigma}^{2k}}\left(\sum
\limits_{i=1}^{n}\left(a_i^{2k}-b_i^{2k}\right)\right)\right\}.
\end{eqnarray*}
It follows that the limiting distribution of $p_n(k)$ is normal if
and only if
\begin{equation}\label{gene1}
\sum \limits_{k=2}^{\infty}
B_{2k}\frac{t^{2k}}{2k(2k)!{\sigma}^{2k}}\left(\sum
\limits_{i=1}^{n}\left(a_i^{2k}-b_i^{2k}\right)\right)\rightarrow 0
\quad \mbox{as} \quad n\rightarrow \infty,
\end{equation}
for $t$ from any bounded set. By virtue of the variance formula
(\ref{var}),  the condition (\ref{gene1}) is equivalent to
\begin{equation}\label{gene2}
\sum \limits_{k=1}^{\infty} B_{2k}\frac{t^{2k}}{2k(2k)!}\frac{\sum
\limits_{i=1}^{n}\left(a_i^{2k}-b_i^{2k}\right)}{\left(\sum
\limits_{i=1}^{n}\left(a_i^2-b_i^2\right)\right)^k}\rightarrow 0
\quad \mbox{as} \quad n\rightarrow \infty
\end{equation}
for $t$ from any bounded set.  Thus (\ref{gene1}) is verified. This
completes the proof.  \qed

\begin{corollary}\label{geco}
Let $p_n(k)$ be given as in the above theorem. Suppose that for
$k\geq 2$, there exist constants $\alpha>0$, $\beta<0$ and
$\gamma<0$ such that
\begin{equation}\label{gene}
\frac{\sum
\limits_{i=1}^{n}\left(a_i^{2k}-b_i^{2k}\right)}{\left(\sum
\limits_{i=1}^{n}(a_i^2-b_i^2)\right)^k}< n^{\gamma}(\alpha
n^{\beta})^{2k},
\end{equation}
for $t$ from any bounded set. Then the limiting distribution of
$p_n(k)$ is normal.
\end{corollary}

 \pf
Note that the convergent radius $R$ of the series
\[ \sum \limits_{2\nmid
k,k\geq3}|B_{2k}|\frac{x^{2k}}{2k(2k)!}
\] is $2\pi$. If  (\ref{gene})
holds for $k>1$, then for $t$ from any bounded set, and  for
sufficiently large $n$, we have
\[
\left|t^{2k}\sum
\limits_{i=1}^{n}\left(a_i^{2k}-b_i^{2k}\right)/{\sigma}^{2k}\right|\leq
n^{\gamma} (t\alpha n^{\beta})^{2k},
\]
where $t\alpha n^{\beta}<2\pi$. It is clear that
$n^{\gamma}\rightarrow 0$ since $\gamma<0$. \qed

If we choose $\alpha=32\sqrt{3}/3$, $2\beta=\gamma=-\frac{1}{3}$,
Theorem \ref{geco} contains  Theorem \ref{main} as a special case.
We now give two more examples. One is the following $q$-analog of
the Catalan numbers
\[
c_n(q)=\frac{[2]}{[2n]}{2n\brack n-1},
\]
which are symmetric and unimodal, see Stanley [].

Using Theorem \ref{gen}, we reach the following assertion.

\begin{corollary}
The distribution of the coefficients in $c_n(q)$ is asymptotically
normal.
\end{corollary}

\noindent \pf  First, we write $c_n(q)$ in the following form:
\[
\frac{\prod\limits_{i=3}^{n}(1-q^{n+i-1})}{(1-q)\prod\limits_{i=3}^{n-1}
(1-q^i)},
\]
Set $a_1=a_2=1,\ a_i=n+i-1,\  3\leq i\leq n$, and $b_1=b_2=1,\
b_3=1,\ b_i=i-1,\ 4\leq i\leq n.$ Then we have
\begin{align*}
\sum_{i=1}^{n}\left(a_i^{2k}-b_i^{2k}\right)=\left(a_3^{2k}-b_3^{2k}\right)+\sum_{i=4}^{n}\left(a_i^{2k}-b_i^{2k}\right)\\[5pt]
=(n+2)^{2k}-1+\sum_{i=3}^{n-1}\left((n+i)^{2k}-i^{2k}\right)
\end{align*}
and
\begin{align*}
\left(\sum_{i=1}^{n}(a_i^2-b_i^2)\right)^k&=\left((n+2)^2-1+\sum_{i=3}^{n-1}\left((n+i)^2-i^2\right)\right)^k\\[5pt]
&=(n-1)^k(n+1)^k(2n-3)^k.
\end{align*}
By the same arguments as in the proof of Lemma \ref{asy1}, we may
set  $\alpha=32\sqrt{3}/3$ and $2\beta=\gamma=-\frac{1}{3}$ such
that the condition (\ref{gene}) is satisfied. Therefore, Theorem
\ref{gen} implies the limiting distribution of the coefficients of
$c_n(q).$\qed

The {\it $m$-Catalan numbers} are defined by
$$C_{n,m}=\frac{1}{(m-1)n+1}{{mn}\choose n},$$ for $n\geq 1$.
Accordingly, the generalized $q$-Catalan numbers are given by
$$C_{n,m}(q)=\frac{1}{[(m-1)n+1]}{{mn}\brack n}.$$
Theorem \ref{gen}  has the following consequence.

\begin{corollary}
The coefficients of the generalized $q$-Catalan numbers $
C_{n,m}(q)$ are normally distributed when $n\rightarrow \infty$.
\end{corollary}
\pf First, express $C_{n,m}(q)$ as follows
\[
\prod \limits_{i=2}^{n}\frac{1-q^{(m-1)n+i}}{1-q^{i}}.
\]
Set $a_1=1,\ a_i=(m-1)n+i,\  2\leq i\leq n$, and $b_1=1,\ \ b_i=i,\
2\leq i\leq n.$ Then we have
\begin{align*}
\sum_{i=1}^{n}\left(a_i^{2k}-b_i^{2k}\right)=\sum_{i=2}^{n}\left(a_i^{2k}-b_i^{2k}\right)
=\sum_{i=2}^{n} \sum \limits_{j=1}^{2k}{{2k}\choose j}
\left((m-1)n\right)^{2k-j}i^{j}.
\end{align*}
The same argument as in the proof of Lemma \ref{asy1} yields the
following bound
\[
\sum_{i=1}^{n}\left(a_i^{2k}-b_i^{2k}\right)<8^{2k}
\left(\left(m-1\right)n\right)^{2k+1}.
\]
Now,
\begin{align*}
\left(\sum_{i=1}^{n}(a_i^2-b_i^2)\right)^k&=\left(\sum_{i=2}^{n}\left(((m-1)n+i)^2-i^2\right)\right)^k\\[5pt]
&>(m-1)^{2k}n^{2k}(n-1)^k\\[5pt]
&>(m-1)^{2k+1}n^{3k}/(2m)^{k}.
\end{align*}
It follows that
\[
\frac{\sum_{i=1}^{n}\left(a_i^{2k}-b_i^{2k}\right)}{\left(\sum_{i=1}^{n}(a_i^2-b_i^2)\right)^k}<(8\sqrt{2m})^{2k}n^{1-k}.
\]

Again, by the same arguments as in the proof of Lemma \ref{asy1}, we
may  set $\alpha=8\sqrt{2m}$ and $2\beta=\gamma=-\frac{1}{3}$ such
that the condition (\ref{gene}) holds. Finally, we may use Theorem
\ref{gen} to get the desired distribution.  \qed

\section{Open Problems}

While the $q$-Catalan numbers are not unimodal for small $n$, see
Stanley \cite{stanlog}, the limiting distribution suggests that the
coefficients are almost unimodal in certain sense for sufficiently
large $n$. Obviously,  the first and the last term should not be
taken into account  otherwise one can never expect to have
unimodality. In fact, an easy computation indicates that $C_n(q)$
are unimodal for $n\geq 16$.

\begin{conjecture}\label{unimodal}
The sequence $\{m_n(1),\ldots ,m_{n}(n(n-1)-1)\}$ is unimodal when
$n$ is sufficiently large.
\end{conjecture}

When $n>70$, numerical evidence  suggestive of  a stronger
conjecture:

\begin{conjecture}\label{log}
There exists an integer $t$ such that  when $n$ is sufficiently
large, the sequence $\{m_n(t),\ldots ,m_n(n(n-1)-t)\}$ is
log-concave, namely,
\[
\left(m_n(k)\right)^2\geq m_n(k+1)m_n(k-1)
\]
for $t+1\leq k\leq n(n-2)-t-1$. Moreover, the minimum value of $t$
seems to be $75$.
\end{conjecture}

We also conjecture that similar properties  hold for the generalized
$q$-Catalan numbers.

\vskip 8pt
 \noindent {\bf Acknowledgments.} We would like to thank
B.H. Margolius for helpful comments. This work was supported by the
973 Project, the PCSIRT Project of the Ministry of Education, the
Ministry of Science and Technology, and the National Science
Foundation of China.

\bibliographystyle{amsplain}

\end{document}